\documentclass[11pt]{article}
\usepackage{amssymb}
\usepackage{amsthm}
\usepackage{latexsym}
\usepackage{graphicx}
\usepackage{fancyhdr}
\usepackage{bbm,bm}
\usepackage{amsmath}
\usepackage{indentfirst}
\usepackage{latexsym,amssymb,amsxtra,mathrsfs,times}
\input xy
\xyoption{all}
\allowdisplaybreaks

\newtheorem{thm}{Theorem}[section]
\newtheorem{lem}[thm]{Lemma}
\newtheorem{cor}[thm]{Corollary}

\theoremstyle{definition}
\newtheorem{defn}{Definition}
\newtheorem{rem}{Remark}

\begin{document}
	\baselineskip=14pt
	
	\title{$\theta$-Congruent Numbers, Tiling Numbers and the Selmer Rank of Related Elliptic Curves: odd $n$}
	
	\author{Qiuyue Liu$^1$, Jing Yang$^{2,}$\footnote{Corresponding author.}\ \ and Keqin Feng$^2$ \\
		{\small $1$ Center for Applied Mathematics, Tianjin University, Tianjin, 300072, China.\ \ \textit{17864309562@163.com.}}\\ 
		{\small $2$ Department of Mathematical Sciences, Tsinghua University, Beijing, 100084, China.}\\
		{\small \textit{y-j@tsinghua.edu.cn,\ \ fengkq@tsinghua.edu.cn.}} }
	\date{}
	
	\maketitle
	
	\begin{abstract}
		Several discrete geometry problems are closely related to the arithmetic theory of elliptic curves defined on the rational fields $\mathbb Q$. In this paper we consider the $\theta$-congruent number for ($\theta=\frac{\pi}{3}$ and $\frac{2\pi}{3}$) and tiling number $n$. For the case that $n\geqslant2$ is square-free odd integer, we determine all $n$ such that the Selmer rank of elliptic curve $E_{n,\frac{\pi}{3}}:\ y^2=x(x-n)(x+3n)$ or/and \(E_{n,\frac{2\pi}{3}}:\ y^2=x(x+n)(x-3n)\) is zero. From this, we provide several series of non $\theta$-congruent numbers for $\theta=\frac{\pi}{3}$ and \(\frac{2\pi}{3}\), and non tiling numbers $n$ with arbitrary many of prime divisors.
		\\[3mm]
		\textbf{Key words}:\ \(\theta\)-congruent number, tiling number, elliptic curve, Selmer rank, odd graph.
		
	\end{abstract}

\section{Introduction}\label{sec-into}

A positive number \(n\) is called a congruent number if it is the area of a right rational triangle, where ``rational'' means that the length of all sides of the triangle are rational numbers. The congruent number problem is to find a simple criterion to determine whether a given integer is a congruent number or not. It is well known that $n$ is congruent if and only if the (Mordell-Weil) rank of the group \(E_n(\mathbb Q)\) of rational points is positive for the following elliptic curve $E_n$ (see Koblitz \cite[Chapter 1]{6}):
\[E_n:\ y^2=x(x+n)(x-n).\]
Many series of congruent and non-congruent numbers has been determined by using arithmetic theory of elliptic curves and modular form theory developed in past half century (see \cite{6,10,11,12} and many others).

Fujiware \cite{3} defined the generalized concept, a $\theta$-congruent number by considering rational triangles with an angle $\theta$. For such a triangle, $\cos\theta$ is rational number: $\cos\theta=\frac{s}{r},\ s,r\in\mathbb{Z},\ \gcd(s,r)=1$ and \(r\geqslant1\). Then $\sin\theta=\frac{\sqrt{r^2-s^2}}{r}$.  

\begin{defn}
 An integer \(n\geqslant1\) is called $\theta$-congruent number if \(n\sqrt{r^2-s^2}\) is the area of a rational triangle with an angle $\theta$.
\end{defn}

For \(\theta=\frac{\pi}{2}\), we have $r=1$ and \(s=0\). Hence $\frac{\pi}{2}$-congruent numbers are the usual congruent numbers. Since $n$ is $\theta$-congruent number if and only if $c^2\theta$ is $\theta$-congruent for some integer $c\neq0$, we may assume without loss of generality that $n$ is square-free. The $\theta$-congruent numbers are also connected with the following elliptic curves:
\[E_{n,\theta}:\ y^2=x(x+(r+s)n)(x-(r-s)n).\]

\begin{thm}[\cite{3}]
  Let $n$ be a square-free positive integer. Then
\par (1). \(n\) is \(\theta\)-congruent if and only if $E_{n,\theta}$ has a rational point of order greater than 2.
\par (2). For \(n\neq1,2,3,6,\ n\) is $\theta$-congruent if and only if the rank of the group $E_{n,\theta}(\mathbb{Q})$ of rational points on $E_{n,\theta}$ is positive.
\end{thm}

In this paper, we deal with the cases $\theta=\frac{\pi}{3}$ and $\frac{2\pi}{3}$. Since $\cos\frac{\pi}{3}=\frac{1}{2}$ and $\cos\frac{2\pi}{3}=-\frac{1}{2}$, we have
\[E_{n,\frac{\pi}{3}}:\ y^2=x(x-n)(x+3n)\ \mbox{ and }\ E_{n,\frac{2\pi}{3}}:\ y^2=x(x+n)(x-3n).\]


Recently, a tiling problem in discrete geometry has been found \cite{7} to be closely related to the elliptic curves $E_{n,\frac{\pi}{3}}$ and $E_{n,\frac{2\pi}{3}}$.

\begin{defn}
	An integer $n\geqslant1$ is called tiling number (TN) if there exists a triangle $\triangle$ and a positive integer $k$ such that $nk^2$ copies of $\triangle$ can be tiled into an equilateral triangle. 
\end{defn}

Following figures show that 1,2,3,4 and 6 are tiling numbers. From the definition we know that for any positive integer $n$ and $m$, $n$ is a tiling number if and only if $nm^2$ is a tiling number. Thus we can assume that $n$  is square-free integer.

\begin{center}
\begin{picture}(40,40)
\put(0,10){\line(1,0){40}}
\put(0,10){\line(1,2){20}}
\put(40,10){\line(-1,2){20}}
\put(5,0){$n=1$}
\end{picture}
\qquad
\begin{picture}(40,40)
\put(0,10){\line(1,0){40}}
\put(0,10){\line(1,2){20}}
\put(40,10){\line(-1,2){20}}
\put(20,10){\line(0,1){40}}
\put(5,0){$n=2$}
\end{picture}
\qquad
\begin{picture}(40,40)
\put(0,10){\line(1,0){40}}
\put(0,10){\line(1,2){20}}
\put(40,10){\line(-1,2){20}}
\put(20,50){\line(0,-1){27}}
\put(0,10){\line(3,2){20}}
\put(40,10){\line(-3,2){20}}
\put(5,0){$n=3$}
\end{picture}
\qquad
\begin{picture}(40,40)
\put(0,10){\line(1,0){40}}
\put(0,10){\line(1,2){20}}
\put(40,10){\line(-1,2){20}}
\put(10,30){\line(1,0){20}}
\put(20,10){\line(1,2){10}}
\put(20,10){\line(-1,2){10}}
\put(5,0){$n=4$}
\end{picture}
\qquad
\begin{picture}(40,40)
\put(0,10){\line(1,0){40}}
\put(0,10){\line(1,2){20}}
\put(40,10){\line(-1,2){20}}
\put(20,50){\line(0,-1){40}}
\put(0,10){\line(3,2){30}}
\put(40,10){\line(-3,2){30}}
\put(5,0){$n=6$}
\end{picture}
\end{center}

\begin{thm}[\cite{7}]\label{thm1.2}
	Let $n\geqslant5$ be a square-free integer. Then the following four statemenets are equivalent to each others.
	\par (1). $n$ is TN;
	\par (2). $n$ is $\frac{\pi}{3}$-CN or $\frac{2\pi}{3}$-CN;
	\par (3). At least one of elliptic curves $E_{n,\frac{\pi}{3}}:y^2=x(x-n)(x+3n)$ and $E_{n,\frac{2\pi}{3}}:y^2=x(x+n)(x-3n)$ have rational solution $(x,y),\ x,y\in\mathbb{Q},y\neq0$.
	\par (4). $\mathrm{rank}\left(E_{n,\frac{\pi}{3}}(\mathbb{Q})\right)\geqslant1$ or $\mathrm{rank}\left(E_{n,\frac{2\pi}{3}}(\mathbb{Q})\right)\geqslant1$.
\end{thm}
\vskip2mm
\noindent
\textbf{Example}. The elliptic curve $E_{5,\frac{2\pi}{3}}:y^2=x(x+5)(x-15)$ has solution $(x,y)=(-1,\pm8)$. By Theorem\ref{thm1.2}, 5 is $\frac{2\pi}{3}$-CN and TN. 

\underline{Conjecture} and \underline{known results}.
\vskip 2mm

Based on deep arithmetic theory on elliptic curves, the following conjecture is raised by Yoshida \cite{10}.
\\
\textbf{Conjecture}. Let $n\geqslant5$ be a square-free integer. Then
\par (1). $n$ is non $\frac{\pi}{3}$-CN\ \ $\Leftrightarrow$\ \ $n\equiv1,2,3,5,7,9,14,15,19\pmod{24}$.
\par (2). $n$ is non $\frac{2\pi}{3}$-CN\ \ $\Leftrightarrow$\ \ $n\equiv1,2,3,6,7,11,13,14,18\pmod{24}$.

\vskip2mm
Several series of $\frac{\pi}{3}$-CN and $\frac{2\pi}{3}$-CN have been found. For all prime numbers $p\equiv23\pmod{24}$, $p$ is $\frac{\pi}{3}$-CN and $\frac{2\pi}{3}$-CN \cite{6}, and $2p,\ 3p$ are $\frac{2\pi}{3}$-CN \cite{13}. Therefore they are TN. On the other hand, for all square-free $n$ with at most two prime divisors $\geqslant5$, Goto \cite{5} computed the Selmer rank of $E_{n,\frac{\pi}{3}}(\mathbb{Q})$ and $E_{n,\frac{2\pi}{3}}(\mathbb{Q})$ listed in \cite{4}, Table 5, 6, 7. It is known that if the Selmer rank of $E_{n,\frac{\pi}{3}}(\mathbb{Q})\ (\mbox{ or }E_{n,\frac{2\pi}{3}}(\mathbb{Q}))$ is zero, then $\mathrm{rank}(E_{n,\frac{\pi}{3}}(\mathbb{Q}))$ (or $E_{n,\frac{2\pi}{3}}(\mathbb{Q})$) is zero and $n$ is non $\frac{\pi}{3}$-CN (non $\frac{\pi}{3}$-CN). With this method Goto obtained many series of non $\frac{\pi}{3}$-CN and/or non $\frac{2\pi}{3}$-CN $n$ for the case that $n$ has at most two prime divisors $\geqslant5$. Goto \cite{5} presented a necessary and sufficient condition for all $n,\ \gcd(n,6)=1$, such that the Selmer rank of $E_{n,\frac{\pi}{3}}(\mathbb{Q})$ is zero, so that $n$ is non $\frac{\pi}{3}$-CN. Such non $\frac{\pi}{3}$-CN $n$ has arbitrary many prime divisors and is described by graph language. (We restate this result as Theorem \ref{thm3.2} in Section \ref{sec3}). 

In this paper, we present a description of odd $n$ in term of graph language such that the Selmer rank of $E_{n,\frac{\pi}{3}}(\mathbb{Q})$ or/and $E_{n,\frac{2\pi}{3}}(\mathbb{Q})$ is zero. Particularly, we give many series of odd non $\frac{\pi}{3}$-CN, non $\frac{2\pi}{3}$-CN and non TN with arbitrary many prime divisors. Comparing with Goto's result in \cite{4} for $\gcd(6,n)=1$ case, we add the case $n\equiv5\pmod{24}$. The results on case $\gcd(n,6)=3$ are new. Besides, we describe all results by an unified-type graph $G(n)$ defined in Section \ref{sec3}. We state and prove our results in Section \ref{sec4} and 5 for case $\gcd(n,6)=1$ and $\gcd(n,6)=3$  respectively. Before proving main results, we introduce some basic facts on elliptic curves and graph theory in Section \ref{sec2} and \ref{sec3}.

\section{Selmer Group of Elliptic Curves}\label{sec2}

In this section, we briefly introduce some basic facts on Selmer group and Selmer rank of elliptic curves. For more detail we refer to \cite{9}.

Consider an elliptic curve
\[E:\ y^2=x^3+Ax^2+Bx\ \ A,B\in\mathbb{Z}\]
with discriminant $\Delta(E)=B^2(A^2-4B)\neq0$. We denote $E(\mathbb{Q})$ the set of rational points $\{(x,y)=(a,b)|a,b\in\mathbb{Q},b^2=a^3+Aa^2+Ba\}$ plus a infinite point $O$. $(E(\mathbb{Q}),\oplus)$ is a finitely generated abelian group with zero element $O$ and a certain operation $\oplus$. Therefore
\[E(\mathbb{Q})=E_t\oplus E_f,\]
where $E_t$ is the finite torsion subgroup of $E(\mathbb{Q})$ and $E_f$ is isomorphic to $(\mathbb{Z}^r,+)$ where $r=\mathrm{rank}(E(\mathbb{Q}))\geqslant0$ is called the rank of $E(\mathbb{Q})$. Therefore, $\mathrm{rank}(E(\mathbb{Q}))=0$ means that $E(\mathbb{Q})$ is a finite group. To determine $\mathrm{rank}(E(\mathbb{Q}))$ is not easy in general, but by using the Selmer group of $E$, we can find some elliptic curves with rank zero.

There is another elliptic curve associated to $E$:
\[E':\ y^2=x^3+A'x^2+B'x,\ \ A'=-2A,\ B'=A^2-4B.\]
The discriminant of $E'$ is $\Delta(E')=16B(A^2-4B)^2\neq0$.

Let
\[M=M(E)=\{\infty\}\cup\{\mbox{prime number: }p:\ p\mid\Delta(E)\}\ (=M(E'));\]
\[D=D(E)=<-1,p:\ p\mid\Delta(E)>\subseteq\mathbb{Q}^*/{\mathbb{Q}^*}^2,\]
where $D(E)$ is the subgroup of $\mathbb{Q}^*/{\mathbb{Q}^*}^2$ generated by $-1$ and all prime divisors of $\Delta(E)$.
For each $p\in M$, we have the following homomorphism of groups
\[\delta_p:E'(\mathbb{Q})\to\mathbb{Q}^*_p/{\mathbb{Q}^*_p}^2,\ \ \ \delta'_p:E(\mathbb{Q})\to\mathbb{Q}_p^*/{\mathbb{Q}_p^*}^2,\]
where
\[\delta_p(P)=\left\lbrace\begin{array}{cl}
1,&\mbox{if } P=O;\\
B',&\mbox{if } P=(0,0);\\
x,&\mbox{if } P=(x,y)\neq(0,0),
\end{array}\right.\ \ \ \ 
\delta'_p(P)=\left\lbrace\begin{array}{cl}
1,&\mbox{if } P=O;\\
B,&\mbox{if } P=(0,0);\\
x,&\mbox{if } P=(x,y)\neq(0,0).
\end{array}\right.\]
And $\mathbb{Q}_p$ is the $p$-adic number fields for prime number $p$, $\mathbb{Q}_\infty=\mathbb{R}$ (the fields of real numbers). We know that
\[\mathbb{Q}^*_p/{\mathbb{Q}^*_p}^2=\left\lbrace\begin{array}{ll}
\{1,g,p,gp\} \mbox{ where } g\in\mathbb{Z},\left(\frac{g}{p}\right)=-1,&\mbox{for prime number } p\geqslant3;\\
\{1,3,5,7,2,6,10,14\},&\mbox{for } p=2;\\
\mathbb{R}^*/{\mathbb{R}^*}^2=\{\pm1\},&\mbox{for } p=\infty.
\end{array}\right.\]
$\mathrm{Im}(\delta_p)$ and $\mathrm{Im}(\delta'_p)$ are subgroups of $\mathbb{Q}^*_p/{\mathbb{Q}^*_p}^2$. From $\mathbb{Q}\subseteq\mathbb{Q}_p$ and $\mathbb{Q}^2\subseteq\mathbb{Q}^2_p$, we can view the elements of $\mathbb{Q}^*/{\mathbb{Q}^*}^2$ as elements in $\mathbb{Q}_p^*/{\mathbb{Q}_p^*}^2$ in natural way.

\begin{defn}
	The Selmer groups of $E(\mathbb{Q})$ are the following two subgroups of $D\ (\subseteq\mathbb{Q}^*/{\mathbb{Q}^*}^2)$:
	\[S(E)=\{d\in D:\ d\in \mathrm{Im}(\delta_p) \mbox{ for each }p\in M\}.\]
	\[S'(E)=\{d\in D:\ d\in \mathrm{Im}(\delta'_p) \mbox{ for each }p\in M\}.\]
\end{defn}

Let $rk_2S(E)$ and $rk_2S'(E)$ be the 2-rank of the elementary 2-groups $S(E)$ and $S'(E)$. It is known that $rk_2S(E)+rk_2S'(E)\geqslant2.$ And
\[\mbox{s-rank}(E)=rk_2S(E)+rk_2S'(E)-2\ (\geqslant 0)\]
is called the Selmer rank of $E(\mathbb{Q})$.

\begin{thm}[\cite{9}]\label{thm2.1}
	$\mathrm{rank}(E(\mathbb{Q}))\leqslant\mbox{s-}\mathrm{rank}(E)$. Particularly, if s-$\mathrm{rank}(E)=0$, then $\mathrm{rank}(E(\mathbb{Q}))=0$.
\end{thm}

$E_{n,\frac{\pi}{3}}$ and $E_{n,\frac{2\pi}{3}}$ belong to the elliptic curves of following type:
\[E=E^{(\alpha,\beta)}:\ y^2=x(x-\alpha)(x-\beta),\]
where $\alpha,\beta\in\mathbb{Z}$ and $0,\alpha,\beta$ are distinct. In his thesis \cite{4}, Goto computed $\mathrm{Im}(\delta_p)$ and $\mathrm{Im}(\delta'_p)$ for $E^{(\alpha,\beta)}$ and then presented an explicit description of $\mathrm{Im}(\delta_p)$ and $\mathrm{Im}(\delta'_p)$ for $E_{n,\frac{\pi}{3}}$ (where $\alpha=n,\beta=-3n$) as following result. 

\begin{thm}[\cite{5}, Lemma 2.8 \& 2.9]\label{thm2.2}
	Let $n\geqslant2$ be a square-free intger. For elliptic curve $E=E_{n,\frac{\pi}{3}},\ p\in M(E)=\{\infty,2,3,p:\ p\mid n\}$ and $d\in D(E)=<-1,2,3,p\mid n>\subseteq\mathbb{Q}^*/{\mathbb{Q}^*}^2.$
\par (A). For $\mathrm{Im}(\delta'_p)$,
  \begin{itemize}
 	\item[(1).] $d\in \mathrm{Im}(\delta'_\infty)$.
 	\item[(2).] $\begin{array}{l}
               \mbox{If } 2\mid d, d\in\mathrm{Im}(\delta'_2)\ \ \Leftrightarrow\ \ 2\mid n \mbox{ and }d\equiv n \pmod8;\\
               \mbox{If } 2\nmid d, d\not\in\mathrm{Im}(\delta'_2)\ \ \Leftrightarrow\ \ n\equiv2,5,6\pmod8 \mbox{ and }d\equiv 3 \pmod4.
            	\end{array}$ 
    \item[(3).] $\begin{array}{l}
    \mbox{If } 3\mid d, d\not\in\mathrm{Im}(\delta'_3)\ \ \Leftrightarrow\ \ n\equiv6\pmod9 \mbox{ and }\frac{d}{3}\equiv1\pmod3;\\
    \mbox{If } 3\nmid d, d\not\in\mathrm{Im}(\delta'_3)\ \ \Leftrightarrow\ \ n\equiv6\pmod9 \mbox{ and }d\equiv2\pmod3.
    \end{array}$ 
    \item[(4).] For the prime divisors $p\geqslant5$ of $n$, \\
    $\begin{array}{l}
    \mbox{If } p\mid d, d\not\in\mathrm{Im}(\delta'_p)\ \ \Leftrightarrow\ \ p\equiv1\pmod3 \mbox{ and }\left(\dfrac{n/d}{p}\right)=-1;\\
    \mbox{If } p\nmid d, d\not\in\mathrm{Im}(\delta'_p)\ \ \Leftrightarrow\ \ p\equiv1\pmod3 \mbox{ and }\left(\dfrac{d}{p}\right)=-1.
    \end{array}$ 
  \end{itemize}
\vskip2mm
\par (B). For $\mathrm{Im}(\delta_p)$, 
   \begin{itemize}
   	\item[(1).] $d\in\mathrm{Im}(\delta_\infty)\ \ \Leftrightarrow\ \ d>0.$
   	\item[(2).]  $\mbox{If } 2\mid d, d\not\in\mathrm{Im}(\delta_2)$;\\
   	$\mbox{If } 2\nmid d, d\in\mathrm{Im}(\delta_2)\ \ \Leftrightarrow\ \ \mbox{one of the following conditions hold}.$
   	\begin{itemize}
   		\item[(a).] $n\equiv5\pmod8$;
   		\item[(b).] $n\equiv1,3,7\pmod8$ and $d\equiv1\pmod4$;
   		\item[(c).] $n\equiv2\pmod8$ and $d\equiv1,7\pmod8$;
   		\item[(d).] $n\equiv6\pmod8$ and $d\equiv1,3\pmod8$.
   	\end{itemize}
   \item[(3).] $\begin{array}{l}
   \mbox{If } 3\mid d, d\in\mathrm{Im}(\delta_3)\ \ \Leftrightarrow\ \ n\equiv6\pmod9 \mbox{ and }\frac{d}{3}\equiv1\pmod3;\\
   \mbox{If } 3\nmid d, d\in\mathrm{Im}(\delta_3)\ \ \Leftrightarrow\ \ d\equiv1\pmod3.
   \end{array}$ 
   \item[(4).] For the prime divisors $p\geqslant5$ of $n$, \\
   $\begin{array}{l}
   \mbox{If } p\mid d, d\in\mathrm{Im}(\delta_p)\ \ \Leftrightarrow\ \ p\equiv1\pmod3 \mbox{ and }\left(\dfrac{-n/d}{p}\right)=1;\\
   \mbox{If } p\nmid d, d\in\mathrm{Im}(\delta_p)\ \ \Leftrightarrow\ \ \left(\dfrac{d}{p}\right)=1.
   \end{array}$ 
   \end{itemize}
 \end{thm}

\vskip2mm\noindent
\textbf{Remark}. As indicated in \cite{4}, there exists dual relationship between $\mathrm{Im}(\delta_p)$ and $\mathrm{Im}(\delta'_p)$ with respect to the Hilbert symbol $(\ ,\ )_p:\mathbb{Q}_p^*/{\mathbb{Q}_p^*}^2\times\mathbb{Q}_p^*/{\mathbb{Q}_p^*}^2\ \to\ \{\pm1\}$ for each $p$ (including $p=\infty$). Therefore $\mathrm{Im}(\delta_p)=\mathrm{Im}(\delta'_p)^\perp$ can be determined from $\mathrm{Im}(\delta'_p)$.
\vskip2mm

With the same way, we can determine $\mathrm{Im}(\delta'_p)$ and $\mathrm{Im}(\delta_p)$ for elliptic curve $E_{n,\frac{2\pi}{3}}$. We state the result and omit the computational details.

\begin{thm}\label{thm2.3}
	Let $n\geqslant2$ be a square-free integer. For elliptic curve $E=E_{n,\frac{2\pi}{3}},\ p\in M(E)=\{\infty,2,3,p: p\mid n\}$ and $d\in D(E)=<-1,2,3,p: p\mid n>\subseteq\mathbb{Q}^*/{\mathbb{Q}^*}^2$, 
\par (A). For $\mathrm{Im}(\delta'_p)$,
   \begin{itemize}
   	\item[(1).] $d\in\mathrm{Im}(\delta'_\infty)$.
   	\item[(2).]  $\begin{array}{l}
   	\mbox{If } 2\mid d, d\in\mathrm{Im}(\delta'_2)\ \ \Leftrightarrow\ \ 2\mid n \mbox{ and }d\equiv-n\pmod8;\\
   	\mbox{If } 2\nmid d, d\not\in\mathrm{Im}(\delta'_2)\ \ \Leftrightarrow\ \ n\equiv2,3,6\pmod8 \mbox{ and }d\equiv 3 \pmod4.
   	\end{array}$ 
   	\item[(3).] $\begin{array}{l}
   	\mbox{If } 3\mid d, d\not\in\mathrm{Im}(\delta'_3)\ \ \Leftrightarrow\ \ n\equiv3\pmod9 \mbox{ and }\frac{d}{3}\equiv1\pmod3;\\
   	\mbox{If } 3\nmid d, d\not\in\mathrm{Im}(\delta'_3)\ \ \Leftrightarrow\ \ n\equiv3\pmod9 \mbox{ and }d\equiv2\pmod3.
   	\end{array}$ 
   	\item[(4).] For the prime divisors $p\geqslant5$ of $n$, \\
   	$\begin{array}{l}
   	\mbox{If } p\mid d, d\not\in\mathrm{Im}(\delta'_p)\ \ \Leftrightarrow\ \ p\equiv1\pmod3 \mbox{ and }\left(\dfrac{-n/d}{p}\right)=-1;\\
   	\mbox{If } p\nmid d, d\not\in\mathrm{Im}(\delta'_p)\ \ \Leftrightarrow\ \ p\equiv1\pmod3 \mbox{ and }\left(\dfrac{d}{p}\right)=-1.
   	\end{array}$ 
   \end{itemize}
\par (B). For $\mathrm{Im}(\delta_p)$,
    \begin{itemize}
   	\item[(1).] $d\in\mathrm{Im}(\delta_\infty)\ \ \Leftrightarrow\ \ d>0.$
   	\item[(2).]  $\mbox{If } 2\mid d, d\not\in\mathrm{Im}(\delta_2);$\\
   	 $\mbox{If } 2\nmid d, d\in\mathrm{Im}(\delta_2)\ \ \Leftrightarrow\ \ \mbox{one of the following conditions hold}.$ 
   	\begin{itemize}
   		\item[(a).] $n\equiv3\pmod8$;
   		\item[(b).] $n\equiv1,5,7\pmod8$ and $d\equiv1\pmod4$;
   		\item[(c).] $n\equiv2\pmod8$ and $d\equiv1,3\pmod8$;
   		\item[(d).] $n\equiv6\pmod8$ and $d\equiv1,7\pmod8$.
   	\end{itemize}
   	\item[(3).] 
   	$\mbox{If } 3\mid d, d\in\mathrm{Im}(\delta_3)\ \ \Leftrightarrow\ \ n\equiv3\pmod9 \mbox{ and }\frac{d}{3}\equiv1\pmod3;$\\
   	$\mbox{If } 3\nmid d, d\in\mathrm{Im}(\delta_3)\ \ \Leftrightarrow\ \ d\equiv1\pmod3.$ 
   	\item[(4).] For the prime divisors $p\geqslant5$ of $n$, \\
   	$\begin{array}{l}
   	\mbox{If } p\mid d, d\in\mathrm{Im}(\delta_p)\ \ \Leftrightarrow\ \ p\equiv1\pmod3 \mbox{ and }\left(\dfrac{n/d}{p}\right)=1;\\
   	\mbox{If } p\nmid d, d\in\mathrm{Im}(\delta_p)\ \ \Leftrightarrow\ \ \left(\dfrac{d}{p}\right)=1.
   	\end{array}$ 
   \end{itemize}
\end{thm}

Theorem \ref{thm2.2} and \ref{thm2.3} can be used to determine the Selmer groups of $E_{n,\theta}(\mathbb{Q})$ for $\theta=\frac{\pi}{3}$ and $\frac{2\pi}{3}$. It is known that $\{1,n,-3n,-3\}\subseteq S'(E_{n,\frac{\pi}{3}})$ and $\{1,-n,3n,-3\}\subseteq S'(E_{n,\frac{2\pi}{3}})$. In section \ref{sec4} and section \ref{sec5} we determine all odd square-free $n$ such that $S(E_{n,\frac{\pi}{3}})=\{1\}$ and $S'(E_{n,\frac{\pi}{3}})=\{1,n,-3n,-3\}$, which implies that s-rank$(E_{n,\frac{\pi}{3}})=0$ and $\mathrm{rank}(E_{n,\frac{\pi}{3}})=0$ by Theorem \ref{thm2.1}. In this way we get a series of non $\frac{\pi}{3}$-CN. Similarly, we determine all odd square-free $n$ such that $S(E_{n,\frac{2\pi}{3}})=\{1\}$ and $S'(E_{n,\frac{2\pi}{3}})=\{1,-n,3n,-3\}$ which implies that $\mathrm{rank}(E_{n,\frac{2\pi}{3}})=0$ and such $n$ is non $\frac{2\pi}{3}$-CN. Then we get a series of non TN.

These non $\frac{\pi}{3}$-CN, non $\frac{2\pi}{3}$-CN and non TN have arbitrary many prime divisors and in terms of graph theory language. In next section, we introduce some basic facts on graph theory we need in last two sections.

\section{Odd Graphs}\label{sec3}

Let $\Gamma=(V,E)$ be a simple directed graph where $V=\{v_1,v_2,\cdots,v_m\}$ is the set of vertices, $E$ is the set of arcs.

A non-trivial partition of $V$ is $\{V_1,V_2\}$ where $V_1\cup V_2=V,\ V_1\cap V_2=\emptyset$ and $V_1, V_2\neq\emptyset$. We denote $V_2=\overline{V}_1$ and view $\{V_1,V_2\}, \{V_2,V_1\}$ as the same partition. There are $2^{m-1}$ partitions of $V$ where $m=|V|$, and one of them is trivial partition.

A partition $\{V_1,V_2\}$ of $V$ is called odd if either there exists $v_1\in V_1$ such that $\#\{v_1\to V_2\}=\#\{v_1\to v:v\in V_2\}$ is odd, or there exists $v_2\in V_2$ such that $\#\{v_2\to V_1\}$ is odd. Otherwise, the partition $\{V_1,V_2\}$ is called even. The trivial partition is even.

Let $\pi=(V_1,\overline{V}_1)$ and $\pi'=(V'_1,\overline{V}'_1)$ be two partitions of $V$. We define a new partition $\pi\vartriangle\pi'=(V^*,\overline{V}^*)$ where
\[V^*=V_1\vartriangle V'_1=(V_1\cup V'_1)\setminus(V_1\cap V'_1).\]

Let $\mathcal{G}$ be the set of all partitions of $V$, $|\mathcal{G}|=2^{m-1}$. Then $(\mathcal{G},\vartriangle)$ is a 2-elementary group with respect to the (Boolean) operation $\vartriangle$. All even partitions of $V$ form a subgroup of $\mathcal{G}$. Trivial partition is the zero element of the group $\mathcal{G}$.

For a simple directed graph $\Gamma=(V,E),\ V=\{v_1,\cdots,v_m\}$, the adjacency matrix is the following $m\times m$ matrix
\[A=A(\Gamma)=(a_{ij})_{1\leqslant i,j\leqslant m}, \ \ \ a_{ij}=\left\lbrace
\begin{array}{cl}
 1,&\mbox{if\ \ }\overrightarrow{v_iv_j}\in E; \\ 0, & \mbox{otherwise.}
\end{array}\right.\]

For each $i\ (1\leqslant i \leqslant m)$, $d_i=\sum_{j=1}^{m}a_{ij}=\#\{v_i\to V\}$ is called the outdegree of $v_i$. The Laplace matrix of $\Gamma$ is defined by 
\[L=L(\Gamma)=\begin{bmatrix}
d_1&&\\&\ddots&\\&&d_m
\end{bmatrix}-A(\Gamma).\]

Since the sum of entries in each row is zero. We get $\mathrm{rank}_{\mathbb{F}_2}L\leqslant m-1$, as we view $L$ to be a matrix over $\mathbb{F}_2$.

\begin{lem}[\cite{1}]\label{lem3.1}
	The number of even partitions of $V$ is $2^{m-l-1}$, where $l=\mathrm{rank}_{\mathbb{F}_2}L.$
\end{lem}

\begin{defn}
	A simple directed graph $\Gamma=(V,E)$ is called odd, if all non-trivial partitions of $V$ are odd. 
\end{defn}

From Lemma \ref{lem3.1}, we know that $\Gamma$ is odd if and only if $\mathrm{rank}_{\mathbb{F}_2}L(\Gamma)=m-1$, where $m=|V|$.

\vskip2mm\noindent
\textbf{Examples}. The following graphs are odd graphs.
	\begin{itemize}
		\item Directed cycle $\overrightarrow{C_m}\ (m\geqslant3):\ V=\{v_1,\cdots,v_m\}, E=\{\overrightarrow{v_1v_2},\overrightarrow{v_2v_3},\cdots, \overrightarrow{v_{m-1}v_m},\overrightarrow{v_{m}v_1}\}.$
		\item Cycle $C_m\ (2\nmid m\geqslant3):\ V=\{v_1,\cdots,v_m\},E=\{\overline{v_1v_2},\overline{v_2v_3},\cdots,\overline{v_{m-1}v_m},\overline{v_mv_1}\}$, where $\overline{v_iv_j}$ means $\overrightarrow{v_iv_j}$ and $\overrightarrow{v_jv_i}$.
		\item Complete graph $K_m\ (2\nmid m\geqslant3):\ V=\{v_1,\cdots,v_m\},E=\{\overline{v_iv_j}:1\leqslant i\neq j\leqslant m\}$.
		\item Trees (connected non-directed graph without cycle).
	\end{itemize}


As early as in 1930's, the conception of odd graphs has been used to determine the 4-rank of the class group of quadratic number fields \cite{10}. A series of non congruent numbers ($\theta=\frac{\pi}{2}$) are presented in terms of odd graphs in \cite{1,2}. With the same method, Goto \cite{5} showed following results on non $\frac{\pi}{3}$-CN.

\begin{thm}[Goto \cite{5}]\label{thm3.2}
	Let $n=p_1\cdots p_t\equiv1,7,19\pmod{24}$ where $p_i\geqslant5\ (1\leqslant i\leqslant t)$ are distinct prime numbers, $t\geqslant1$. We define two graphs:
	\[G(-3n):=(V_1,E_1),\ V_1=\{-1,3,p_1,\cdots,p_t\},\ E_1=\left\{\overrightarrow{p_iv}:v\in V, \left(\frac{v}{p_i}\right)=-1, p_i\neq v,1\leqslant i\leqslant t\right\};\]
	\[g(n):=(V_2,E_2),\ V_2=\{-1,p_1,\cdots,p_t\},\ E_2=\left\{\overrightarrow{p_ip_j}:\left(\frac{p_j}{p_i}\right)=-1,1\leqslant i\neq j\leqslant t \right\}\bigcup\left\{\overline{p_i(-1)}:\left(\frac{-1}{p_i}\right)=-1,1\leqslant i\leqslant t\right\}.\]
	Then for $E=E_{n,\frac{\pi}{3}}$, 
 \par (1). $S'(E)=\{1,n,-3n,-3\}$ if and only if
   \par\quad(a). $p_i\equiv1\pmod3\ (1\leqslant i\leqslant t)$, and 
   \par\quad(b). In the graph $G(-3n)$, all partitions $(V_1,\overline{V}_1)$ of $V=\{-1,3,p_1\cdots,p_t\}$ are odd except the trivial partition and $\pi=\{V_1=\{-1,3\},\overline{V}_1=\{p_1,\cdots,p_t\}\}$. 
 \par (2). Assume that $p_i\equiv1\pmod3\ (1\leqslant i\leqslant t)$, then $S(E)=\{1\}$ if and only if $g(n)$ is odd graph.
\end{thm}

In next two sections we determine all odd square-free integers $n$ such that
  \begin{itemize}
	\item[(a).] The Selmer rank s-$\mathrm{rank}(E_{n,\frac{\pi}{3}})$ is zero, so that $n$ is non $\frac{\pi}{3}$-CN. or,
	\item[(b).] The Selmer rank s-$\mathrm{rank}(E_{n,\frac{2\pi}{3}})$ is zero, so that $n$ is non $\frac{2\pi}{3}$-CN. 
\end{itemize}
From (a) and (b) we present a series of odd non tiling numbers. We describe such odd square-free $n$ by using the following an unified-type graphs.

Let $m=p_1\cdots p_t$ or $m=-p_1\cdots p_t\ (t\geqslant1)$ where $p_1,\cdots,p_t$ are distinct prime numbers. We define the graph $G(m)=(V,E)$ by
\[V=\left\lbrace\begin{array}{ll}
  \{p_1,\cdots,p_t\},&\mbox{if } m>0;\\   \{-1,p_1,\cdots,p_t\},&\mbox{if } m<0.
   \end{array}\right.\ \ \mbox{ and }\ \ 
   E=\left\{\overrightarrow{p_iv}:v\in V,\left(\frac{v}{p_i}\right)=-1,p_i\equiv1\pmod3,1\leqslant i\leqslant t\right\}.\]
   
For a nontrivial partition $\{V_1,V_2=\overline{V}_1\}$ of $V$, let
\[d=d(V_1)=\prod_{v\in V_1}v,\ \ \ \bar{d}=\prod_{v\in V_2}v.\]
Then $m=\prod_{v\in V}v=d\bar{d}$. For $v\in V_1,\ v\equiv1\pmod3$, we have $\#\{v\to V_2\}=\#\left\{v'\in V_2:\left(\frac{v'}{v}\right)=-1\right\}$. Therefore $\#\{v\to V_2\}$ is odd $\Leftrightarrow$ $\prod_{v'\in V_2}\left(\frac{v'}{v}\right)=\left(\frac{\bar{d}}{v}\right)$ is odd. And then, $\{V_1, V_2=\overline{V}_1\}$ is odd partition $\Leftrightarrow$ there exists a prime number $v\mid m,\ v\equiv1\pmod3$ such that either $v\mid d$ (means $v\in V_1$) and $\left(\frac{m/d}{v}\right)=-1$, or $v\nmid d$ (means $v\in V_2$) and $\left(\frac{d}{v}\right)=-1.$

This basic fact is made possible to describe $n$ being non $\frac{\pi}{3}$-CN or non $\frac{2\pi}{3}$-CN by oddness of the graph $G(m)$ with certain $m\in\{\pm n,\pm3n\}$.

\section{The case $\bm{\gcd(n,6)=1}$}\label{sec4}

\begin{thm}\label{thm4.1}
	Let $n=p_1\cdots p_t\ (t\geqslant1)$ where $p_i\geqslant5\ (1\leqslant i\leqslant t)$ are distinct prime numbers. $E=E_{n,\frac{\pi}{3}}$.
	\begin{itemize}
		\item[(I).] For $n\equiv1,7,19\pmod{24}$.
		\begin{itemize}
			\item[(I.1)] $S'(E)=\{1,n,-3n,-3\}$ if and only if the following two conditions hold.
			\begin{itemize}
				\item[(a).] $p_i\equiv1\pmod3\ (1\leqslant i\leqslant t)$;
				\item[(b).] In graph $G(-3n)$, all partitions of $V=\{-1,3,p_1,\cdots,p_t\}$ are odd except trivial partition and $\pi=\{V_1=\{-1,3\},V_2=\{p_1,\cdots,p_t\}\}$. 
			\end{itemize} 
		  \item[(I.2)] Assume that $p_i\equiv1\pmod3\ (1\leqslant i\leqslant t)$. Then $S(E)=\{1\}$ if and only if in graph $G(-n)$, for all $V_1\subseteq\{p_1,\cdots,p_t\}$ satisfying $d(V_1)>1$ and $d(V_1)\equiv1\pmod4$, the partitions $\{V_1,\overline{V}_1\}$ of $V=\{-1,p_1,\cdots,p_t\}$ is odd.
		\end{itemize} 
	  \item[(II).] For $n\equiv5\pmod{24}$,
	  \begin{itemize}
	  	\item[(II.1)] $S'(E)=\{1,n,-3n,-3\}$ if and only if the following two conditions hold.
	  	\begin{itemize}
	  		\item[(a).] $p_i\equiv1\pmod3\ (1\leqslant i\leqslant t-1)$ and $p_t\equiv2\pmod3$;
	  		\item[(b).] In graph $G(-3n)$, for all $V_1\subseteq V$ satisfying $d(V_1)=1\pmod4$, the partitions of $\{V_1,\overline{V}_1\}$ of $V=\{-1,3,p_1,\cdots,p_t\}$ is odd except trivial partition and $V_1=\{-1,3\}$.
	  	\end{itemize}
  	  \item[(II.2)] Assume that the condition (a) in (II.1) are satisfied. Then $S(E)=\{1\}$ if and only if in graph $G(-n)$, for all $V_1\subseteq\{p_1,\cdots,p_t\},\ V_1\neq\emptyset$, the partition $\{V_1,\overline{V}_1\}$ of $V=\{-1,p_1,\cdots,p_t\}$ is odd. 	  	
	  \end{itemize}
	\end{itemize}
\end{thm}
\noindent
\textbf{\underline{Proof}.} \underline{(I.1)}: This is Theorem \ref{thm3.2} (1).

\underline{(I.2)}: By Theorem \ref{thm3.2} (2), we need to show that the condition on graph $G(-n)$ in (I.2) is equivalent to that the graph $g(n)$ is odd. Recall that the vertices set of both $g(n)$ and $G(-n)$ is $V=\{-1,p_1,\cdots,p_t\}$. For a nontrivial partition $\pi=\{V_1,V_2=\overline{V}_1\}$ of $V$, we can assume $-1\in V_2$. It is easy to see that if $d(V_1)\equiv3\pmod4$, then in graph $g(n)$, the number $\#\{-1\to V_1\}$ is odd, and then the partition $\pi$ is odd. If $d(V_1)\equiv1\pmod4$, then the number $\#\{-1\to V_1\}$ is even in $g(n)$. Thus the partition $\pi$ is odd in $g(n)$ if and only if $\pi$ is odd in $G(-n)$. This shows that the condition of $G(-n)$ in (I.2) is equivalent that $g(n)$ is odd graph.

\underline{(II.1)} Assume that $n=p_1\cdots p_rq_1\cdots q_s\equiv5\pmod{24},\ r+s=t\geqslant1$, and $p_i\equiv1\pmod3\ (1\leqslant i\leqslant r)$ and $q_j\equiv2\pmod3\ (1\leqslant j \leqslant s)$ are distinct odd prime numbers. From $n\equiv2\pmod3$, we know that $t=r+s$ is odd.

In order $S'(E_{n,\frac{\pi}{3}})=\{1,n,-3n,-3\}$, we need to show that for each $d\in<-1,2,3,p_1,\cdots,p_t>\setminus\{1,n,-3n,-3\}$ where exists $p\in\{\infty,2,3,p_1,\cdots,p_t\}$ such that $d\not\in \mathrm{Im}(\delta'_p)$.

By Theorem \ref{thm2.2} (A), we know that $d\in \mathrm{Im}(\delta'_\infty),\ d\in\mathrm{Im}(\delta'_3)$. And by $n\equiv5\pmod{24}$, we get $d\in\mathrm{Im}(\delta'_2)\ \Leftrightarrow\ d\equiv1\pmod4$. Therefore we need for each $d\equiv1\pmod4$, there exists $p\in\{p_1,\cdots,p_r,q_1,\cdots,q_s\}$ such that $d\not\in\mathrm{Im}(\delta'_p)$. Namely (by Theorem \ref{thm2.3} (A.4)), $p\equiv1\pmod3$, and $\left(\frac{n/d}{p}\right)=-1$ when $p\mid d$ or $\left(\frac{d}{p}\right)=-1$ when $p\nmid d$. As we explained at the end of last section, this just means that in graph $G(-3n)$, all partitions $\{V_1,\overline{V}_1\}$ of $V=\{-1,3,p_1,\cdots,p_r,q_1,\cdots,q_s\}$ stisfying $d(V_1)\equiv1\pmod4$ are odd except $V_1=\emptyset$ and $V_1=\{-1,3\}$. (Remark: from $p_i\equiv1\pmod3\ (1\leqslant i\leqslant r)$, we get $\left(\frac{-1}{p_i}\right)\left(\frac{3}{p_i}\right)=\left(\frac{-3}{p_i}\right)=1$. This implies that in graph $G(-3n)$, the partition $\{V_1=\{-1,3\},\overline{V}_1\}$ is even.)

Moreover, $d(V_1)d(\overline{V}_1)=-3n\equiv1\pmod4$, we get $d(V_1)\equiv d(\overline{V}_1)\pmod4$. Thus the set 
\[T=\{(V_1,\overline{V}_1):d(V_1)\equiv1\pmod4\}\]
is a subgroup of $(\mathcal{G},\vartriangle)$, where $\mathcal{G}$ is the group of the all partitions of $V=\{-1,3,p_1,\cdots,p_r,q_1,\cdots,q_s\}$. $\mathcal{G}$ has two coset $T$ and $\pi\vartriangle T$ where $\pi=\{V_1=\{-1,3\},\overline{V}_1\}$. Since $T$ has two even partitions, the coset $\pi\vartriangle T$ has at most two even partitions. Therefore the number of even partitions in $\mathcal{G}$ is at most four. By Lemma \ref{lem3.1}, for the $(t+2)\times(t+2)$ Laplace matrix $L$ of the graph $G(-3n),\ \mathrm{rank}_{\mathbb{F}_2}L\geqslant(t+2)-3=t-1=r+s-1$. On the other hand, there is no arc $\overrightarrow{vv'}$ in $G(-3n)$ for $v\in\{-1,3,q_1,\cdots,q_s\}$. This means that there are $s+2$ rows in $L$ being zero vectors. Thus $r+s-1\leqslant\mathrm{rank}_{\mathbb{F}_2}L\leqslant(r+s-2)-(s+2)=r$. We get $s=1$, $n=p_1\cdots,p_t,\ p_i\equiv1\pmod3\ (1\leqslant i\leqslant t-1),\ p_t\equiv2\pmod3$ and the graph $G(-3n)$ satisfies the condition (II.1)(b).

\underline{(II.2)} Assume that $n=p_1\cdots p_t\equiv5\pmod{24},\ p_i\equiv1\pmod3\ (1\leqslant i \leqslant t-1)$ and $p_t\equiv2\pmod3$. We want to determine $n$ such that $S(E_{n,\frac{\pi}{3}})=\{1\}$. Namely, we need to show that for each $d\in<-1,2,3,p_1,\cdots,p_t>\subseteq\mathbb{Q}^*/{\mathbb{Q}^*}^2,\ d\neq1$, there exists $p\in\{\infty,2,3,p_1\cdots,p_t\}$ such that $d\not\in\mathrm{Im}(\delta_p)$.

By Theorem \ref{thm2.2}(B), we know that $d<0\ \Rightarrow\ d\not\in\mathrm{Im}(\delta_\infty);\ 2\mid d\ \Rightarrow\ d\not\in\mathrm{Im}(\delta_2);\ 3\mid d\ \Rightarrow d\not\in\mathrm{Im}(\delta_3)$. 
Thus we consider $d\in<p_1,\cdots,p_t>$ and $d\neq1$. 
In this case, $d\in\mathrm{Im}(\delta_2)$, and $d\in\mathrm{Im}(\delta_3)\ \Leftrightarrow\ d\equiv1\pmod3$. 
Therefore, for $d\equiv1\pmod3$ and $d>1$, we need that there exists $p\in\{p_1,\cdots,p_{t-1}\}$ such that $d\not\in\mathrm{Im}(\delta_p)$. 
By Theorem \ref{thm2.2}(B.4), this means that either $p\mid d$ and $\left(\frac{-n/d}{p}\right)=-1$, or $p\in\{p_1,\cdots,p_{t-1},p_t\}$, $p\nmid d$ and $\left(\frac{d}{p}\right)=-1$. For the last case, if $\left(\frac{d}{p_t}\right)=-1$ we may take $p=p_t$, otherwise we need to have $p\in\{p_1,\cdots,p_{t-1}\}$ such that $p\nmid d$ and $\left(\frac{d}{p}\right)=-1$. Now we consider the graph $G(-n)$ with $V=\{-1,p_1,\cdots,p_t\}$. It is easy to see that 
\[T=\left\{(V_1,\overline{V}_1):V_1\subseteq\{p_1,\cdots,p_{t-1}\},\left(\frac{d(V_1)}{p_t}\right)=1\right\}\] 
is a subgroup of $(\mathcal{G},\vartriangle)$ and the above requirements is just that in the graph $G(-n)$, all nontrivial partitions in $T$ are odd. 

This completes the proof of Theorem \ref{thm4.1}. \hfill$\square$

\begin{cor}\label{cor4.2}
	Let $n=p_1\cdots,p_t\ (t\geqslant1)$ where $p_i\geqslant5\ (1\leqslant i\leqslant t)$ are distinct prime numbers, $E=E_{n,\frac{\pi}{3}}$. Then
	\begin{itemize}
		\item[(I).] For $n\equiv1,7,19\pmod{24}$, the Selmer rank of $E(\mathbb{Q})$ is zero, so that $n$ is non $\frac{\pi}{3}$-CN if and only if $p_i\equiv1\pmod3\ (1\leqslant i\leqslant t)$ and the graphs $G(-3n)$ and $G(-n)$ satisfy the condition (I.1)(b) and (I.2) of Theorem \ref{thm4.1} respectively.
		\item[(II).] For $n\equiv5\pmod{24}$, the Selmer rank of $E(\mathbb{Q})$ is zero, so that $n$ is non $\frac{\pi}{3}$-CN if and only if $p_i\equiv1\pmod3\ (1\leqslant i\leqslant t-1),\ p_t\equiv2\pmod3$ and the graphs $G(-3n),G(-n)$ satisfy the condition (II.1)(b) and (II.2) of Theorem \ref{thm4.1} respectively. \hfill$\square$
	\end{itemize}
\end{cor}

The following result provides non $\frac{\pi}{3}$-CN $n$ with arbitrary many of prime divisors.

\begin{cor}\label{cor4.3}
	(I). Let $n=p_1\cdots p_t\equiv7,19\pmod{24}$, where $p_i\equiv1\pmod3\ (1\leqslant i\leqslant t)$ are distinct prime numbers. Then the Selmer rank of $E_{n,\frac{\pi}{3}}$ is zero if and only if the graph $G(n)$ is odd.
	\par (II). Let $n=p_1\cdots p_t\equiv5\pmod{24}$, where $p_i\equiv1\pmod{12}\ (1\leqslant i\leqslant t-1)$ and $p_t\equiv5\pmod{12}$ are distinct prime numbers. If the graph $G(n)$ is odd, then $n$ is non $\frac{\pi}{3}$-CN.
\end{cor} 
\noindent
\underline{\textbf{Proof.}} (I). Suppose that $G(n)$ is an odd graph, we need to show that $G(-3n)$ and $G(-n)$ satisfy the condition (I.1)(b) and (I.2) of Theorem \ref{thm4.1} respectively.

For the graph $G(-3n)$ where $V=\{-1,3,p_1,\cdots,p_t\}$, let $\pi=\{V_1,V_2\}$ be a partition of $V$. If $\{-1,3\}\subseteq V_1$, then $V_1=\{-1,3\}\cup V'_1$ and $\pi'=\{V'_1,V'_2\}$ is a partition of the vertices set $V'=\{p_1,\cdots,p_t\}$ of $G(n)$. $\pi$ is nontrivial and $\pi\neq\{\{-1,3\},\{p_1,\cdots,p_t\}\}$ if and only if $\pi'$ is nontrivial. Moreover, from $p_i\equiv1\pmod3$ we get $1=\left(\frac{-3}{p_i}\right)$, therefore $\left(\frac{-1}{p_i}\right)=\left(\frac{3}{p_i}\right)\ (1\leqslant i\leqslant t)$. Then this implies that partition $\pi$ is odd in $G(-3n)$ if and only if partition $\pi'$ is odd in $G(n)$. By assumption that $G(n)$ is odd graph, we know that partition $\pi$ is odd. Next, we assume that $-1\in V_1$ and $3\in V_2$. Then $V_1=\{-1\}\cup V'_1,\ V_2=\{3\}\cup V'_2$, where $\pi'=\{V'_1,V'_2\}$ is a partition of $V'=\{p_1,\cdots,p_t\}$. Now we show that such partition $\pi=\{V_1,V_2\}$ is odd in graph $G(-3n)$. For dong this, consider the total number $\Sigma(\pi)$ of arcs between $V_1$ and $V_2$.
\[\begin{array}{rcl}
  \Sigma(\pi)&=&\#\{V_1\to V_2\}+\#\{V_2\to V_1\}\\
   &=&\#\{V'_1\to\{3\}\}+\#\{V'_1\to V'_2\}+\#\{V'_2\to\{-1\}\}+\#\{V'_2\to V'_1\}.
\end{array}\]

Let $P=\{1\leqslant i\leqslant t: p_i\equiv1\pmod4\},\ Q=\{1\leqslant i\leqslant t: p_i\equiv3\pmod4\},\ |P|=r,|Q|=s$. Then $r+s=t, 2\nmid s$ since $n\equiv3\pmod4$. Let $Q_i=Q\cap V'_i,\ |Q_i|=s_i\ (i=1,2)$. Then $s_1+s_2=s\equiv1\pmod2$ so that $s_1s_2\equiv0\pmod2$. For each $p\in\{p_1,\cdots,p_t\}=P\cup Q$, from $p\equiv1\pmod3$ we get $\left(\frac{-3}{p}\right)=1$ so that
\[\left(\frac{3}{p}\right)=\left(\frac{-1}{p}\right)=\left\lbrace 
  \begin{array}{ll}
   1,&\mbox{if }p\in P;\\ -1,&\mbox{if }p\in Q.
  \end{array}\right. \]
Therefore $\#\{V'_1\to\{3\}\}=|Q_1|=s_1,\ \#\{V'_2\to \{-1\}\}=|Q_2|=s_2$. By quadratic reciprocity law, for $p\in P\cap V'_1$, then $\#\{p\to V'_2\}=\#\{V'_2\to p\}$, for $p\in P\cap V'_2$, $\#\{p\to V'_1\}=\#\{V'_1\to p\}$, and for $p,p'\in Q,\ \#\{p\to p'\}+\#\{p'\to p\}=1$. Therefore
\[\Sigma(\pi)\equiv s_1+s_2+\#\{Q_1\to Q_2\}+\#\{Q_2\to Q_1\}\equiv s_1+s_2+s_1s_2\equiv1\pmod2.\]
From this we know that the partition $\pi$ is odd. Thus the graph $G(-3n)$ satisfies the condition (I.1)(b) of Theorem \ref{thm4.1}.

Next we show that $G(n)$ is odd if and only if $G(-n)$ is odd. Let $n=p_1\cdots p_rq_1\cdots q_s\ (r+s=t)$ where $p_i\equiv1\pmod4\ (1\leqslant i\leqslant r)$ and $q_j\equiv3\pmod4\ (1\leqslant j\leqslant s)$. From $n\equiv3\pmod4$ we get $2\nmid s$. Suppose that the adjacency and Laplace matrix of $G(n)$ are (they are viewed as $t\times t$ matrices over $\mathbb{F}_2$) 
\[A(n)=\begin{array}{cc}
 & \begin{array}{cccccc} p_1&\cdots&p_r&q_1&\cdots&q_s\end{array} \\
 \begin{array}{c} p_1\\\vdots\\p_r\\q_1\\\vdots\\q_s\end{array} & 
 \left[ \begin{array}{cccccc} &&&&&\\ &&&&&\\ 
 \multicolumn{3}{c}{\raisebox{2ex}[0pt]{\Huge$\ \ P'$\ }} & \multicolumn{3}{c}{\raisebox{2ex}[0pt]{\Huge\ $A$\ \ }}\\
 	&&&&&\\ &&&&&\\[3mm] 
 	\multicolumn{3}{c}{\raisebox{2ex}[0pt]{\Huge $\ \ A\ $}} & \multicolumn{3}{c}{\raisebox{2ex}[0pt]{\Huge\ $Q'$\ \ }}
 		\end{array}\right] 
\end{array},
\ \ L(n)= \left[\begin{array}{cccccc} &&&&&\\ &&&&&\\ 
\multicolumn{3}{c}{\raisebox{2ex}[0pt]{\Huge$\ P$\ }} & \multicolumn{3}{c}{\raisebox{2ex}[0pt]{\Huge\ $A$\ }}\\
&&&&&\\ &&&&&\\[3mm] 
\multicolumn{3}{c}{\raisebox{2ex}[0pt]{\Huge $\ A\ $}} & \multicolumn{3}{c}{\raisebox{2ex}[0pt]{\Huge\ $Q$\ }}
\end{array}\right] .\]
Then the adjacency and Laplace matrices of $G(-n)$ over $\mathbb{F}_2$ are
\[A(-n)=\begin{array}{rc}
& \begin{array}{ccccccc} -1&p_1&\cdots&p_r&q_1&\cdots&q_s\end{array} \\
\begin{array}{c} -1\\p_1\\\vdots\\p_r\\q_1\\\vdots\\q_s\end{array} & 
\left[\begin{array}{c|ccc|ccc} 
0&0&\cdots&0&0&\ \cdots&0\\\hline
0&&&&&&\\ \vdots&&&&&&\\ 
0&\multicolumn{3}{c|}{\raisebox{2ex}[0pt]{\Huge$P'$}} & \multicolumn{3}{c}{\raisebox{2ex}[0pt]{\Huge$A$}}\\\hline
1&&&&&&\\ \vdots&&&&&&\\
1&\multicolumn{3}{c|}{\raisebox{2ex}[0pt]{\Huge $A$}} & \multicolumn{3}{c}{\raisebox{2ex}[0pt]{\Huge$Q'$}}
\end{array}\right] 
\end{array},\ \ 
L(-n)=\left[\begin{array}{c|ccc|ccc} 
0&0&\cdots&0&0&\quad\cdots&0\\\hline
0&&&&&&\\ \vdots&&&&&&\\ 
0&\multicolumn{3}{c|}{\raisebox{2ex}[0pt]{\Huge$P$}} & \multicolumn{3}{c}{\raisebox{2ex}[0pt]{\Huge$A$}}\\\hline
1&&&&&&\\ \vdots&&&&&&\\
1&\multicolumn{3}{c|}{\raisebox{2ex}[0pt]{\Huge $A$}} & \multicolumn{3}{c}{\raisebox{2ex}[0pt]{\Huge$Q+I_s$}}
\end{array}\right].\]
where $I_s$ is the $s\times s$ identity matrix. By quadratic reciprocity law, we get $P=P^T,\ A=A^T$ and $Q+Q^T=I_s+J_s$, where $J_s$ is the $s\times s$ matrix with all entries being 1. Thus the transposition of $L(-n)$ is 
\[L(-n)^T=\left[\begin{array}{c|ccc|ccc} 
0&0&\cdots&0&1&\quad\cdots&1\\\hline
0&&&&&&\\ \vdots&&&&&&\\ 
0&\multicolumn{3}{c|}{\raisebox{2ex}[0pt]{\Huge$P$}} & \multicolumn{3}{c}{\raisebox{2ex}[0pt]{\Huge$A$}}\\\hline
0&&&&&&\\ \vdots&&&&&&\\
0&\multicolumn{3}{c|}{\raisebox{2ex}[0pt]{\Huge $A$}} & \multicolumn{3}{c}{\raisebox{2ex}[0pt]{\Huge$Q+J_s$}}
\end{array}\right].\]
In the Laplace matrix $L(n)$, the summation of each row is zero. Then by $2\nmid s$, add all columns of $L(-n)^T$ to the first one, we get
\[M=\left[\begin{array}{c|ccc|ccc} 
1&0&\cdots&0&1&\quad\cdots&1\\\hline
0&&&&&&\\ \vdots&&&&&&\\ 
0&\multicolumn{3}{c|}{\raisebox{2ex}[0pt]{\Huge$P$}} & \multicolumn{3}{c}{\raisebox{2ex}[0pt]{\Huge$A$}}\\\hline
1&&&&&&\\ \vdots&&&&&&\\
1&\multicolumn{3}{c|}{\raisebox{2ex}[0pt]{\Huge $A$}} & \multicolumn{3}{c}{\raisebox{2ex}[0pt]{\Huge$Q+J_s$}}
\end{array}\right].\]
Then add the first column to each of the last $s$ columns, we get
\[N=\left[\begin{array}{c|ccc|ccc} 
1&0&\cdots&0&0&\cdots&0\\\hline
0&&&&&&\\ \vdots&&&&&&\\ 
0&\multicolumn{3}{c|}{\raisebox{2ex}[0pt]{\Huge$P$}} & \multicolumn{3}{c}{\raisebox{2ex}[0pt]{\Huge$A$}}\\\hline
1&&&&&&\\ \vdots&&&&&&\\
1&\multicolumn{3}{c|}{\raisebox{2ex}[0pt]{\Huge $A$}} & \multicolumn{3}{c}{\raisebox{2ex}[0pt]{\Huge$Q$}}
\end{array}\right].\]
Therefore $\mathrm{rank}_{\mathbb{F}_2}L(-n)=\mathrm{rank}_{\mathbb{F}_2}N=1+\mathrm{rank}_{\mathbb{F}_2}L(n)$, and 
\[\begin{array}{ll}
G(n)\mbox{ is odd graph}\quad&\Leftrightarrow\quad \mathrm{rank}_{\mathbb{F}_2}L(n)=r+s-1\\
&\Leftrightarrow\quad\mathrm{rank}_{\mathbb{F}_2}L(-n)=r+s \\
&\Leftrightarrow\quad G(-n)\mbox{ is odd graph}.\end{array}\]
Particularly, if $G(n)$ is odd graph, then $G(-n)$ satisfies the condition (I.2) of Theorem \ref{thm4.1}.

(II). Let $n=p_1\cdots p_t\equiv5\pmod{24}\ (t\geqslant1)$, $p_i\equiv1\pmod{12}\ (1\leqslant i\leqslant t-1),\ p_t\equiv5\pmod{12}.$ Suppose that $G(n)$ is odd graph. We need to show that $G(-3n)$ and $G(-n)$ satisfy the condition (II.1)(b) and (II.2) of Theorem \ref{thm4.1} respectively.

Let $\pi=\{V_1,V_2\}$ be a partition of $V=\{-1,3,p_1,\cdots,p_t\},\ d(V_1)\equiv1\pmod4,\ V_1\neq\emptyset,\{-1,3\}$. From $d(V_1)\equiv1\pmod4$ and $p_i\equiv1\pmod4\ (1\leqslant i\leqslant t)$ we have that $\{-1,3\}\cap V_1=\emptyset$ or $\{-1,3\}\subseteq V_1$. Then for $V'=\{p_1,\cdots,p_t\},\ V'_1=V_1\cap V'$ and $V'_2=V_2\cap V',\ \pi'=(V'_1,V'_2)$ is a nontrivial partition of $V'$. From assuption that $G(n)$ is odd, we know that the partition $\pi'$ is odd, then by $\left(\frac{-1}{p_i}\right)\left(\frac{3}{p_i}\right)=\left(\frac{-3}{p_i}\right)=1\ (1\leqslant i\leqslant t)$ we know that the partition $\pi$ is odd. This means that the graph $G(-3n)$ satisfies the condition (II.1)(b) of Theorem \ref{thm4.1}. 

Let $\pi=(V_1,V_2)$ be a partition of $V=\{-1,p_1,\cdots,p_t\},\ -1\in V_2,\ V_1\subseteq\{p_1,\cdots,p_t\}$ and $V_1\neq\emptyset$. Then $\pi'=(V_1,V_2\setminus\{-1\})$ is a nontrivial partition of $\{p_1,\cdots,p_t\}$. By assumption that $G(n)$ is odd graph,  we know that the partition $\pi'$ of $\{p_1,\cdots,p_t\}$ is odd. Then the partition $\pi$ is odd since $\left(\frac{-1}{p_i}\right)=1$ for all $i,\ (1\leqslant i\leqslant t)$. This means that graph $G(n)$ satisfies the condition (II.2) of Theorem \ref{thm4.1}.\hfill$\square$
\vskip2mm

Now we consider the Selmer groups of elliptic curve $E_{n,\frac{2\pi}{3}}$.

\begin{thm}\label{thm4.4}
	Let $n=p_1\cdots  p_tq_1\cdots q_s\equiv1,7,11,13\pmod{24}$ where $p_i\equiv1\pmod3\ (1\leqslant i\leqslant t)$ and $q_j\equiv2\pmod3\ (1\leqslant j\leqslant s)$ are $t+s\ (\geqslant1)$ distinct odd prime numbers. Then
	\begin{itemize}
		\item[(I).] $S'(E_{n,\frac{2\pi}{3}})=\{1,-n,3n,-3\}$ if and only if the following two conditions hold
		\begin{itemize}
			\item[(a).] $s=0,\ n=p_1\cdots p_t\equiv1,7,13\pmod{24},\ t\geqslant1;$
			\item[(b).] Graph $G(-n)$ is odd. 
		\end{itemize} 
	   \item[(II).] Let $n=p_1\cdots p_t\equiv1,7,13\pmod{24},\ t\geqslant1$. Then $S(E_{n,\frac{2\pi}{3}})=\{1\}$ if and only if 
	   \begin{itemize}
	   	\item[(a).] $n\equiv7\pmod{24}$; and
	   	\item[(b).] Graph $G(n)$ is odd. 
	   \end{itemize}
	\end{itemize}
\end{thm}
\noindent
\underline{\textbf{Proof.}} (I). By Theorem \ref{thm2.3}(A), for $d\in<-1,2,3,p_1,\dots, p_t,q_1,\dots,q_s>\subseteq\mathbb{Q}^*/{\mathbb{Q}^*}^2$, we have $d\in\mathrm{Im}(\delta'_\infty)$ and $d\in \mathrm{Im}(\delta'_3)$. And $d\in\mathrm{Im}(\delta'_2)\ \Leftrightarrow\ 2\nmid d$. Now we consider $d\in<-1,3,p_1,\dots,p_t,q_1,\dots,q_s>\setminus\{1,-n,3n,-3\}$. In this case, $d\not\in S'(E_{n,\frac{2\pi}{3}})$ if and only if there exists $p\in\{p_1,\dots,p_t\}$ such that either $p\mid d$ and $\left(\frac{-n/d}{p}\right)=-1$, or $p\nmid d$ and $\left(\frac{d}{p}\right)=-1$. It is easy to see that if $d$ satisfies this requirement, so does $-3d$ (since $\left(\frac{-3}{p}\right)=1$). Therefore the above requirement just means that the graph $G(-n)$ is odd (trivial partition of $V=\{-1,p_1,\dots,p_t,q_1,\dots,q_s\}$ corresponds to $d=1,3n,-3$ and $-n$). Moreover, let $L$ be the Laplace matrix of $G(-n)$, then $\mathrm{rank}_{\mathbb{F}_2}L=(t+s+1)-1=t+s$. Since there are $s+1$ rows in $L$ being zero vectors. We get $\mathrm{rank}_{\mathbb{F}_2}L\leqslant(t+s+1)-(s+1)=t$. Therefore $s=0$ and $n=p_1\cdots p_t\ (t\geqslant1)$. At last, from $n=p_1\cdots p_t\equiv1\pmod3$ we get $n\not\equiv11\pmod{24}$.

(II). Suppose that $n=p_1\cdots p_t\equiv1,7,13\pmod{24},\ t\geqslant1$ and $p_i\equiv1\pmod3\ (1\leqslant i\leqslant t)$ are distinct odd prime numbers. In order to get $S(E_{n,\frac{2\pi}{3}})=\{1\}$, we need to find the condition of $n$ such that for each $d\in<-1,2,3,p_1,\dots,p_t>,\ d\neq1$, there exists $p\in\{\infty,2,3,p_1,\dots,p_t\}$ such that $d\not\in\mathrm{Im}(\delta_p)$. From Theorem \ref{thm2.3}(B) we know that $d\in\mathrm{Im}(\delta_\infty)\ \Leftrightarrow\ d>0; 2\mid d \Rightarrow\ d\not\in\mathrm{Im}(\delta_2)$; and $3\mid d\ \Rightarrow\ d\not\in\mathrm{Im}(\delta_3)$. Then we consider $d\in<p_1,\dots,p_t>$. In this case, $d\in\mathrm{Im}(\delta_2)\ \Leftrightarrow\ d\equiv1\pmod4$; $d\in\mathrm{Im}(\delta_3)\ \Leftrightarrow\ d\equiv1\pmod3$ which is always true since $p_i\equiv1\pmod3\ (1\leqslant i\leqslant t)$. Thus we consider $d\mid p_1\cdots p_t,\ d>1$ and $d\equiv1\pmod4$. In this case we need to have $p\in\{p_1,\dots,p_t\}$ such that $d\not\in\mathrm{Im}(\delta_p)$. From Theorem \ref{thm2.3}(B.4), this requirement is just either $p\mid d$ and $\left(\frac{n/d}{p}\right)=-1$, or $p\nmid d$ and $\left(\frac{d}{p}\right)=-1$.

When $n\equiv1,13\pmod{24}$, we have that $n\equiv1\pmod3$. For each $p_i\mid n,\ \left(\frac{n/n}{p_i}\right)=1$ which means $n\in S(E_{n,\frac{2\pi}{3}})$. Thus if $S(E_{n,\frac{2\pi}{3}})=\{1\}$ then $n\equiv7\pmod{24}$. And the above requirement is just that $n\not\in S(E_{n,\frac{2\pi}{3}})$ and in graph $G(n)$ each nontrivial partition $(V_1,V_2)$ of $V=\{p_1,\cdots,p_t\}$ satisfying $d(V_1)\equiv1\pmod4$ is odd. For $n\equiv7\pmod  {24}$ we get $n\equiv3\pmod4$ and $n\not\in S_{E,\frac{2\pi}{3}}$. Moreover, for any partition $(V_1,V_2)$ of $V$, $d(V_1)\cdot d(V_2)=n\equiv3\pmod4$. Therefore $d(V_1)\equiv1\pmod4$ or $d(V_2)\equiv1\pmod4$. Therefore the above requirement just means that $G(n)$ is odd graph. 
\hfill$\square$

\begin{cor}\label{cor4.5}
	Let $n=p_1\cdots p_t\equiv7\pmod{24}$ where $p_i\equiv1\pmod3\ (1\leqslant i\leqslant t, t\geqslant1)$ are distinct odd prime numbers. Then the following statements are equivalent to each others.
	\begin{itemize}
		\item[(1).] $G(n)$ is an odd graph;
		\item[(2).] The Selmer rank of $E_{n,\frac{\pi}{3}}$ is zero (which implies that $n$ is non $\frac{\pi}{3}$-CN);
		\item[(3).] The Selmer rank of $E_{n,\frac{2\pi}{3}}$ is zero (which implies that $n$ is non $\frac{2\pi}{3}$-CN).
		\end{itemize}
	Therefore, if $G(n)$ is an odd graph then $n$ is non TN.
\end{cor}
\noindent
\underline{\textbf{Proof.}} (1)$\Leftrightarrow$(2): this is Corollary \ref{cor4.3}(1).
\par (1)$\Leftrightarrow$(3): By Theorem \ref{thm4.4} and the fact showed in Corollary \ref{cor4.3}(1): $G(n)$ is odd if and only if $G(-n)$ is odd. \hfill$\square$

\begin{rem}
If $n=p_1\cdots p_t\equiv7\pmod{24},\ p_i\equiv1\pmod3\ (1\leqslant i\leqslant t)$ and $G(n)$ is an odd graph, J. Pan and Y. Tian \cite{14} proved that the BSD conjecture for $E_{n,\frac{\pi}{3}}$ and $E_{n,\frac{2\pi}{3}}$ is true. 
\end{rem}

\section{The case \bm{$\gcd(n,6)=3$}}\label{sec5}

\begin{thm}\label{thm5.1}
	Let $n=3p_1\cdots p_tq_1\cdots q_s\equiv3,9,15\pmod{24}$ where $p_i\equiv1\pmod3\ (1\leqslant i\leqslant t),\ q_j\equiv2\pmod3\ (1\leqslant j\leqslant s)$ are distinct odd prime numbers, $t+s\geqslant1$. Let $P=p_1\cdots p_t,\ Q=q_1\cdots q_s,\ E=E_{n,\frac{\pi}{3}}$.
	\begin{itemize}
		\item[(I).] For $S'(E)$.
		\begin{itemize}
			\item[(I.1)] If $PQ\equiv1\pmod3$, then $S'(E)=\{1,n,-\frac{n}{3},-3\}$ if and only if $s=0$ and $G(n)\ (n=3p_1\cdots p_t,t\geqslant1)$ is an odd graph.
			\item[(I.2)] IF $PQ\equiv2\pmod3$, then $S'(E)=\{1,n,-\frac{n}{3},-3\}$ if and only if $s=1$ and for graph $G(n)$, all nontrivial partitions $(V_1,\overline{V}_1)$ of $V=\{3,p_1,\cdots,p_t,q_1\}$ satisfying $V_1\subseteq\{p_1,\cdots,p_t\}$ are odd.
		\end{itemize} 
	\item[(II).] For $S(E)$. 
	\begin{itemize}
		\item[(II.1)] If $PQ\equiv1\pmod3$ and $s=0$ so that $n=3p_1\cdots p_t\ (t\geqslant1)$, then $S(E)=\{1\}$ if and only if $n\equiv9\pmod{24}$ and $G(n/3)$ is an odd graph. 
		\item[(II.2)] If $PQ\equiv2\pmod3$ and $s=1$ so that $n=3p_1\cdots p_tq_1$, then $S(E)=\{1\}$ if and only if for graph $G(-n)$, all non trivial partitions of $V=\{-1,3,p_1,\cdots,p_t,q_1\}$ satisfying $V_1\subseteq\{3,p_1,\cdots,p_t\},\ d(V_1)\equiv1\pmod4$ and $\left(\frac{d(V_1)}{q_1}\right)=1$ are odd. 
	\end{itemize}
	\end{itemize}
\end{thm}
\noindent
\textbf{\underline{Proof}.} (I). $S'(E)=\{1,n,-\frac{n}{3},-3\}$ requires that for each $d\in<-1,2,3,p_1,\cdots,p_t,q_1,\cdots,q_s>,\ d\not\in\{1,n,-\frac{n}{3},-3\}$, there exists $p\in\{\infty,2,3,p_1,\cdots,p_s,q_1,\cdots,q_s\}$ such that $d\not\in\mathrm{Im}(\delta'_p)$. By Theorem \ref{thm2.2}(A), we know that $d\in\mathrm{Im}(\delta'_\infty);\ d\in\mathrm{Im}(\delta'_2)\ \Leftrightarrow\ 2\nmid d$; for $1\leqslant j\leqslant s,\ d\in\mathrm{Im}(\delta'_{q_j})$. Thus it is reduced to consider $d\in<-1,3,p_1,\cdots,p_t,$ $q_1,\cdots,q_s>\setminus\{1,n,-\frac{n}{3},-3\}$ and $p\in\{3,p_1,\cdots,p_s\}$.

(I.1) Assume that $PQ\equiv1\pmod3$. By Theorem \ref{thm2.2}(A), $d\in\mathrm{Im}(\delta'_3)\ \Leftrightarrow\ 3\nmid d$. Thus we consider $d\in<-1,p_1,\cdots,p_t,q_1,\cdots,q_s>\setminus\{1,-\frac{n}{3}\}$ and $p\in\{p_1,\cdots,p_t\}$. In this case, 
$$\begin{array}{rcl}
d\not\in S'(E)&\Leftrightarrow& \mbox{there exists }p\in\{p_1,\cdots,p_t\} \mbox{ such that }d\not\in\mathrm{Im}(\delta'_p)\\
&\Leftrightarrow&\mbox{ there exists }p\in\{p_1,\cdots,p_t\} \mbox{ such that }\left(\frac{n/d}{p}\right)=-1\ (\mbox{for }p|d)\ \mbox{or }\left(\frac{d}{p}\right)=-1\ (\mbox{for }p\nmid d)\end{array}$$
Therefore, $S'(E)=\{1,n,-n/3,-3\}\ \Leftrightarrow\ G(n)$ is an odd graph. (Remark that for $p\equiv1\pmod3$, $\left(\frac{-1}{p}\right)=\left(\frac{3}{p}\right)$).

Moreover, let $L$ be Laplace matrix of $G(n)$ with size $(t+s+1)\times(t+s+1)$. If $G(n)$ is odd, then $\mathrm{rank}_{\mathbb{F}_2}L=t+s$. On the other hand, $L$ has $s+1$ zero-rows since there exists no arc starting from vertices $3,q_1,\cdots,q_s$. We get $t+s=\mathrm{rank}_{\mathbb{F}_2}L\leqslant(t+s+1)-(s+1)=t$ which inplies $s=0$.

(I.2) Assume that $PQ\equiv p_1\cdots p_tq_1\cdots q_s\equiv2\pmod3$. Then $2\nmid s$. By Theorem \ref{thm2.2}(A), 
\[d\in\mathrm{Im}(\delta'_3)\ \Leftrightarrow\ d\equiv1\pmod3\ (\mbox{for }3\nmid d)\mbox{ or }\frac{d}{3}\equiv2\pmod3\ (\mbox{for }3\mid d).\]
Since $-3\in S'(E)$, we know that $d\in S'(E)\ \Leftrightarrow\ -3d\in S'(E)$. Thus we can consider $3\nmid d$,. Namely, $d\in<-1,p_1,\cdots,p_t,q_1,\cdots,q_s>\setminus\{1,-\frac{n}{3}\}$, and $d\equiv1\pmod3$. Then by Theorenm \ref{thm2.2}(A),
\[\begin{array}{ll}
&d\not\in S'(E)\Leftrightarrow\mbox{there exists }p\in\{p_1,\cdots,p_t\}\mbox{ such that } d\not\in\mathrm{Im}(\delta'_p)\\
\Leftrightarrow&\mbox{there exists }p\in\{p_1,\cdots,p_t\} \mbox{ such that }\left(\frac{n/d}{p}\right)=-1\mbox{ if $p|d$ or }\left(\frac{d}{p}\right)=-1\mbox{ if }p\nmid d.  
\end{array}\]
Therefore,
$S'(E)=\{1,n,-\frac{n}{3},-3\}\ \Leftrightarrow\ $For graph $G(-\frac{n}{3})$, all nontrivial partitions $(V_1,\overline{V}_1)$ of $V=\{-1,p_1,\cdots,p_t$, $q_1,\cdots,q_s\}$ satisfying $d(V_1)\equiv1\pmod3$ are odd.

All nontrivial partitions $(V_1,\overline{V}_1)$ with $d(V_1)\equiv1\pmod3$ (and then $d(\overline{V}_1)=\frac{-PQ}{d(V_1)}\equiv1\pmod3$) form a subgroup $T$ of the group $\mathcal{G}$ of all partitions of $V$, and $[\mathcal{G}:T]=2$. In subgroup $T$, only trivial partition is even, thus there exists at most one even partition in another coset. Therefore the total number of even partitions in $\mathcal{G}$ is at most two. This implies that $\mathrm{rank}_{\mathbb{F}_2}L\geqslant t+s-1$ where $L$ is the Laplace matrix of $G(-n/3)$ with size $(t+s+1)\times(t+s+1)$. On the other hand, $L$ has $s+1$ zero-rows, we get $t+s-1\leqslant\mathrm{rank}_{\mathbb{F}_2}L\leqslant(t+s+1)-(s+1)=t$. Therefore $s=1$.

(II). For $S(E)$, we assume that $n=3p_1\cdots p_tq_1\cdots q_s\equiv3,9,15\pmod{24},\ s=0$ or 1. By Theorem \ref{thm2.2}(B), $d\in\mathrm{Im}(\delta_\infty)\ \Leftrightarrow\ d>0$; for $2\mid d,\ d\not\in\mathrm{Im}(\delta_2)$ and for $2\nmid d$, $d\in\mathrm{Im}(\delta_2)\ \Leftrightarrow\ d\equiv1\pmod4$. Thus it is reduced to consider $p\in\{3,p_1,\cdots,p_t,q_1,\cdots,q_s\},\ d\in<3,p_1,\cdots,p_t,q_1,\cdots,q_s>$, and $d\equiv1\pmod4$. 

(II.1) Assume that $PQ\equiv1\pmod3$. Then $s=0$ and $n=3p_1\cdots p_t\ (t\geqslant1)$. By Theorem \ref{thm2.2}(B), $d\in\mathrm{Im}(\delta_3)\ \Leftrightarrow\ 3\nmid d$. Thus we consider $d\in<p_1,\cdots,p_t>,\ d\equiv1\pmod4$ and $p\in\{p_1,\cdots,p_t\}$. In this case,
\[\begin{array}{rcl}
d\not\in S(E)&\Leftrightarrow& \mbox{there exists }p\in\{p_1,\cdots,p_t\}\mbox{ such that }d\not\in\mathrm{Im}(\delta_p)\\
&\Leftrightarrow&\mbox{there exists }p\in\{p_1,\cdots,p_t\}\mbox{ such that }\left(\frac{d}{p}\right)=-1\mbox{ (for $p\nmid d$) or }\left( \frac{-n/d}{p}\right)=-1\mbox{ (for $p\mid  d$)}.  
\end{array}\]
Therefore
\(S(E)=\{1\}\ \Leftrightarrow\) for graph $G(-n)$, all nontrivial partitions $(V_1,\overline{V}_1)$ of $V=\{-1,3,p_1,\cdots,p_t\}$ satisfying $V_1\subseteq\{p_1,\cdots,p_t\},\ d(V_1)\equiv1\pmod4$ are odd. Moreover, if $n=3P\equiv3,15\pmod{24}$, then $P=p_1\cdots p_t\equiv1\pmod4$. The partition ($V_1=\{p_1,\cdots,p_t\},\overline{V}_1=\{-1,3\}$) is even since $p_i\equiv1\pmod3$ and $\left(\frac{-1}{p_i}\right)\cdot\left(\frac{3}{p_i}\right)=\left(\frac{-3}{p_i}\right)=1\ (1\leqslant i\leqslant t)$. Therefore, if $S(E)=\{1\}$, then $n\equiv9\pmod{24}$. In this case, for any  partition $(V_1,\overline{V}_1)$ of $V=\{-1,3,p_1,\cdots,p_t\},\ d(V_1)d(\overline{V}_1)=-3P\equiv3\pmod4$. We know that $d(V_1)\equiv1\pmod4$ or $d(\overline{V}_1)\equiv1\pmod4$. Then the result can be restated as follows.

For $n=3p_1\cdots p_t\equiv9\pmod{24},\ S(E)=\{1\}$ if and only if $G(n/3)$ is an odd graph.

(II.2) Assume that $PQ\equiv2\pmod3$, then $n=3p_1\cdots p_tq_1\equiv3,9,15\pmod{24}$. By Theorem \ref{thm2.2}(B) we know that
\[
d\in\mathrm{Im}(\delta_3)\ \Leftrightarrow\ \frac{d}{3}\equiv1\pmod3\ (\mbox{when }3|d) \mbox{ or }d\equiv1\pmod3\ (\mbox{when }3\nmid d).\]
Then for $d\in<3,p_1,\cdots,p_t>,\ d\neq1,\ d\equiv1\pmod4$, we have that
\[\begin{array}{rcl}
d\not\in S(E)&\Leftrightarrow& \mbox{there exists }p\in\{p_1,\cdots,p_t,q_1\}\mbox{ such that }d\not\in\mathrm{Im}(\delta_p)\\
&\Leftrightarrow& \mbox{there exists }p\in\{p_1,\cdots,p_t\}\mbox{ such that }p\mid d,\ \left(\frac{-n/d}{p}\right)=-1\\
&& \mbox{or there exists }p\in\{p_1,\cdots,p_t,q_1\}\mbox{ such that }p\nmid d,\ \left(\frac{d}{p}\right)=-1.\\
&\Leftrightarrow& \mbox{In graph }G(-n)\mbox{, all nontrivial partitions }(V_1,\overline{V}_1)\mbox{ of }V=\{-1,3,p_1,\cdots,p_t,q_1\}\\
&&\mbox{satisfying }V_1\subseteq\{3,p_1,\cdots,p_t\},\ d(V_1)\equiv1\pmod4\mbox{ and }\left(\frac{d(V_1)}{q_1}\right)=1\mbox{ are odd.} 
\end{array}\]
\hfill $\square$

\begin{thm}\label{thm5.2}
	Let $n=3p_1\cdots p_tq_1\cdots q_s\equiv3\pmod{24}$ where $p_i\equiv1\pmod3\ (1\leqslant i\leqslant t)$ and $q_j\equiv2\pmod3\ (1\leqslant j\leqslant s)$ are distinct prime numbers $\geqslant5$, $t+s\geqslant1$. Let $P=p_1\cdots p_t, Q=q_1\cdots q_s,E=E_{n,\frac{2\pi}{3}}(\mathbb{Q})$.
	\begin{itemize}
		\item[(I).] For $S'(E)$
		\begin{itemize}
			\item[(I.1)] If $PQ\equiv1\pmod3$, then $S'(E)\{1,-n,n/3,-3\}$ if and only if $s=0$ or 2, and when $s=0$, for $G(-n)$, all nontrivial partitions of $V=\{-1,3,p_1,\cdots,p_t\}$ satisfying $d(V_1)\neq-3,P$ and $d(V_1)\equiv1\pmod4$ are odd. When $s=2$, for $G(-n)$, all nontrivial partitions $(V_1,\overline{V}_1)$ of $V=\{-1,3,p_1,\cdots,p_t,q_1,q_2\}$ satisfying $d(V_1)\neq Pq_1q_2,d(V_1)\equiv1\pmod{4}$, and $d(V_1)\equiv1\pmod{12}$ are odd.
			\item[(I.2)] If $PQ\equiv2\pmod3$, then $S'(E)=\{1,-n,n/3,-3\}$ if and only if $s=1$ and for graph $G(-n)$, all nontrivial partitions $(V_1,\overline{V}_1)$ of $V=\{-1,3,p_1,\cdots,p_t,q_1\}$ satisfying $d(V_1)\neq-3,Pq_1$ are odd.
		\end{itemize}
	\item[(II).] For $S(E)$
	\begin{itemize}
		\item[(II.1)] If $PQ\equiv1\pmod3$, and $s\in\{0,2\}$, then $S(E)=\{1\}$ if and only if $s=2$ and for graph $G(n)$, all nontrivial partitions $(V_1,\overline{V}_1)$ of $V=\{3,p_1,\cdots,p_t,q_1,q_2\}$ satisfying $V_1\subseteq\{3,p_1,\cdots,p_t\}$ and $\left(\frac{d(V_1)}{q_1}\right)=\left(\frac{d(V_1)}{q_2}\right)=1$ are odd.
		\item[(II.2)] If $PQ\equiv2\pmod3$, and $s=1$, then $S(E)=\{1\}$ if and only if for graph $G(n)$, all nontrivial partitions $(V_1,\overline{V}_1)$ of $V=\{3,p_1,\cdots,p_t,q_1\}$ satisfying $V_1\subseteq\{p_1,\cdots,p_t\}$ and $\left(\frac{d(V_1)}{q_1}\right)=1$ are odd.
	\end{itemize}
	\end{itemize}
\end{thm}
\noindent
\textbf{\underline{Proof}.} (I.1) Assume that $PQ\equiv1\pmod3$, then $s$ is even. $S'(E)=\{1,-n,n/3,-3\}$ requires that for $d\in<-1,2,3,p_1,\cdots,p_t,q_1,\cdots,q_s>\setminus\{1,-n,n/3,-3\}$  there exists $p\in\{\infty,2,3,p_1,\cdots,p_t,q_1,\cdots,q_s\}$ such that $d\not\in\mathrm{Im}(\delta'_p)$. From Theorem \ref{thm2.3}(A) we know that $d\in\mathrm{Im}(\delta'_\infty)$; $2\mid d\ \Rightarrow\ d\not\in\mathrm{Im}(\delta'_2)$; if $2\nmid d$, then $d\in\mathrm{Im}(\delta'_2)\ \Leftrightarrow\ d\equiv1\pmod4.$ And $d\in\mathrm{Im}(\delta'_3)\ \Leftrightarrow\ \frac{d}{3}\equiv2\pmod3$ (for $3|d$) or $d\equiv1\pmod3$ (for $d\nmid d$).

Since $-3\in S'(E)$, it is reduced to consider $d\in<-1,p_1,\cdots,p_t,q_1,\cdots,q_s>,\ d\neq1,n/3$ and $d\equiv1\pmod{12}$. For this case, by Theorenm \ref{thm2.2}(A) again,
\[\begin{array}{rcl}
d\not\in S'(E)&\Leftrightarrow&\mbox{there exists }p\in\{p_1,\cdots,p_t,q_1,\cdots,q_s\}\mbox{ such that }d\not\in\mathrm{Im}(\delta'_p)\\
 &\Leftrightarrow&\mbox{there exists }p\in\{p_1,\cdots,p_t\}\mbox{ such that } \left(\frac{-n/d}{p}\right)=-1\mbox{ (for $p|d$) or }\left(\frac{d}{p}\right)=-1\mbox{ (for $p\nmid d$)}\end{array}\]
Therefore,
$S'(E)=\{1,-n,n/3,-3\}\ \Leftrightarrow\ \mbox{for graph }G(-n)\mbox{, all nontrivial partitions }(V_1,\overline{V}_1)\mbox{ of }V=\{-1,3,p_1,\cdots,p_t,$ $q_1,\cdots,q_s\}\mbox{ satisfying }d(V_1)\equiv1\pmod{12},\ d(V_1)\neq PQ\mbox{ are odd.}$

Moreover, all partitions $(V_1,\overline{V}_1)$ satisfying $d(V_1)\equiv1\pmod{12}$ form a subgroup of $\mathcal{G}$ with $\varphi(12)=4$ cosets. Each coset has at most two even patitions, so the total number of even partitions in $\mathcal{G}$ is at most $16=2^4$. This implies that $\mathrm{rank}_{\mathbb{F}_2}L\geqslant(2+s+t)-4-1=s+t-3$. On the other hand, the Laplace matrix $L$ of $G(-n)$ has $2+s$ zero-rows which implies that $\mathrm{rank}_{\mathbb{F}_2}L\leqslant(2+t+s)-(2+s)=t$. Therefore $s=0$ or 2. In these two cases, it is to see that graph $G(n)$ holds the requirement stated in (I.1).

(I.2) If $PQ\equiv2\pmod3$, $s$ is odd. From Theorenm \ref{thm2.3}(A) we know that 
\[d\in\mathrm{Im}(\delta'_\infty);\ d\in\mathrm{Im}(\delta'_3);\ d\in\mathrm{Im}(\delta'_2)\ \Leftrightarrow\ 2\nmid d\mbox{ and }d\equiv1\pmod4.\]
Now we consider $d\in<-1,3,p_1,\cdots,p_t,q_1,\cdots,q_s>\setminus\{1,-n,n/3,-3\}$ and $d\equiv1\pmod4$. In this case,
\[d\not\in S'(E)\ \Leftrightarrow\ \mbox{there exists }p\in\{p_1,\cdots,p_t\}\mbox{ such that }\left( \frac{-n/d}{p}\right)=-1\mbox{ (for $p|d$) or }\left(\frac{d}{p}\right)=1\mbox{ (for $p\nmid d$)}.\]
Therefore,
\[\begin{array}{rcl}
S'(E)=\{1,-n,n/3,-3\}&\Leftrightarrow&\mbox{for graph }G(-n)\mbox{, all nontrivial partitions }(V_1,\overline{V}_1)\mbox{ of }\\
&&V=\{-1,3,p_1,\cdots,p_t,q_1,\cdots,q_s\}\mbox{ satisfying }d(V_1)\neq-3,PQ\\
&&\mbox{and }d(V_1)\equiv1\pmod4\mbox{ are odd.}
\end{array}\]
Moreover, all partitions $(V_1,\overline{V}_1)$ satisfying $d(V_1)\equiv1\pmod4$ form a subgroup of $\mathcal{G}$ with two cosets. Each coset has at most two even partitions. The total number of even partitions in $\mathcal{G}$ is at most $4=2^2$. Therefore 
$\mathrm{rank}_{\mathbb{F}_2}L(G(-n))\geqslant(2+t+s)-3=t+s-1$. On the other hand, $L$ has $s+2$ zero-rows which implies that
$\mathrm{rank}_{\mathbb{F}_2}L(G(-n))\geqslant(2+t+s)-(2+s)=t$. Therefore $s=1$ and graph $G(-n)$ holds the requirement stated in (I.2).

(II). For $S(E)$, from Theorem \ref{thm2.3}(B) we know that
\[d\in\mathrm{Im}(\delta_\infty)\ \Leftrightarrow\ d>0;\ d\in\mathrm{Im}(\delta_2)\ \Leftrightarrow\ 2\mid d.\]
Thus it is reduced to consider $d\in<-1,3,p_1,\cdots,p_t,q_1,\cdots,q_s>$ and $p\in\{3,p_1,\cdots,p_t,q_1,\cdots,q_s\}$.

(II.1) Assume that $PQ\equiv1\pmod3$ and $s\in\{0,2\}$. By Theorem \ref{thm2.3}(B), 
\[d\in\mathrm{Im}(\delta_3)\ \Leftrightarrow\ \frac{d}{3}\equiv1\pmod3\mbox{ (when $3|d$) or }d\equiv1\pmod3\mbox{ (when $3\nmid d$)}.\]
If the righ-hand side holds, then 
\[\begin{array}{rcl}
 d\not\in S(E)&\Leftrightarrow&\mbox{there exists }p\in\{p_1,\cdots,p_t,q_1\cdots,q_s\}\mbox{ such that }d\not\in\mathrm{Im}(\delta_p)\\
 &\Leftrightarrow& \mbox{either there exists }p\in\{q_1,\cdots,q_s\}\mbox{ such that $p\nmid d$ and }\left(\frac{d}{p}\right)=-1\\
 &&\mbox{or there exists }p\in\{p_1,\cdots,p_t\}\mbox{ such that } \left(\frac{n/d}{p}\right)=-1\mbox{ or }\left(\frac{d}{p}\right)=-1.
\end{array}\]
If $s=0$, then $n=3p_1\cdots p_t\equiv3\pmod{24}$ and $P=p_1\cdots p_t\equiv1\pmod8$. From above argument we get $P\in S(E)$. Therefore $s=2$. In this case, if $q_1\mid d$ or $q_2\mid d$, then $d\not\in S(E)$. Thus we consider $d\in<3,p_1,\cdots,p_t>$, and 
\[\begin{array}{rcl}
d\not\in S(E)&\Leftrightarrow&\left(\frac{d}{q_1}\right)=-1\mbox{ or }\left(\frac{d}{q_2}\right)=-1\mbox{, or there exists }p\in\{p_1,\cdots,p_t\}\\
&&\mbox{such that }\left(\frac{n/d}{p}\right)=-1\mbox{ or }\left(\frac{d}{p}\right)=-1.
\end{array}\]
Therefore,
\[\begin{array}{rcl}
S(E)=\{1\}&\Leftrightarrow&\mbox{for graph $G(n)$, all nontrivial partitions }(V_1,\overline{V}_1)\mbox{ of }V=\{3,p_1,\cdots,p_t,q_1,q_2\}\\
&&\mbox{satisfying }V_1\subseteq\{3,p_1,\cdots,p_t\},\ \left(\frac{d(V_1)}{q_1}\right)=\left(\frac{d(V_1)}{q_2}\right)=1\mbox{ are odd.}
\end{array}\]

(II.2) Assume that $PQ\equiv2\pmod3$ and $s=1$. Then $n=3PQ\equiv6\pmod9$. By Theorem \ref{thm2.3}(B), $d\in\mathrm{Im}(\delta_3)\ \Leftrightarrow\ 3\nmid d$ and $d\equiv1\pmod3$. With the similar argument we get
\[\begin{array}{rcl}
S(E)=\{1\}&\Leftrightarrow&\mbox{for graph $G(n)$, all nontrivial partitions }(V_1,\overline{V}_1)\mbox{ of }V=\{3,p_1,\cdots,p_t,q_1\}\\
&&\mbox{satisfying }V_1\subseteq\{p_1,\cdots,p_t\}\mbox{ and} \left(\frac{d(V_1)}{q_1}\right)=1\mbox{ are odd.}
\end{array}\]
This completes the proof of Theorem \ref{thm5.2}.\hfill$\square$
\vskip2mm

At the end of this paper, we present a series $n$ of non TN with arbitrary many of prime divisors and $\gcd(n,6)=3$.

\begin{cor}\label{cor5.3}
	Let $n=3p_1\cdots p_tq_1\equiv3\pmod{24}$ where $p_i\equiv1\pmod{12}\ (1\leqslant i\leqslant t)$ and $q_1\equiv5\pmod{12}$ are distinct prime numbers. If $G(\frac{n}{3})$ is an odd graph, then the Selmer rank of $E_{n,\frac{\pi}{3}}$ and $E_{n,\frac{2\pi}{3}}$ is zero. Therefore such $n$ is non $\frac{\pi}{3}$-CN, non $\frac{2\pi}{3}$-CN and non TN.
\end{cor} 
\noindent
\textbf{\underline{Proof}.} By assumption we know that $Pq_1=p_1\cdots p_tq_1\equiv2\pmod3$ and $n\equiv3\pmod{24}$. From $p_i\equiv1\pmod{12}$ we get $\left(\frac{3}{p_i}\right)=\left(\frac{-1}{p_i}\right)=1\ (1\leqslant i\leqslant t)$. If $G(\frac{n}{3})=G(p_1\cdots p_t)$ is an odd graph, it is easy to see that graph $G(n)$ satisfies the requirement in Theorem \ref{thm5.1}(I.2) and Theorem \ref{thm5.2}(II.2) and graph $G(-n)$ satisfies the requirement in Theorem \ref{thm5.1}(II.2) and Theorem \ref{thm5.2}(I.2). Therefore the Selmer rank of both $E_{n,\frac{\pi}{3}}$ and $E_{n,\frac{2\pi}{3}}$ is zero and $n$ is non $\frac{\pi}{3}$-CN, $\frac{2\pi}{3}$-CN and non TN.\hfill$\square$


\begin{thebibliography}{HD}
	
	
	
	
	\normalsize
	\baselineskip=11pt
	
	
	\bibitem{1} K. Feng, Non-Congruent numbers, odd graphs and the BSD conjecture, Acta Arith. 80(1996), 71-83.
	\bibitem{2} K. Feng and M. Xiong, On elliptic curves $y^2=x^3-n^2x$ with rank zero, Jour. Number Theory, 109(2004), 1-26.
	\bibitem{3} M. Fujiwara, $\theta$-congruent numbers, in Number Theory, de Gruyter Berlin, 1998, 235-241.
	\bibitem{4} T. Goto, A study on the Selmer groups of elliptic curves with a rational 2-torsion, Doc. Thesis, Kyushu Univ. 2002, 54pp.
	\bibitem{5} T. Goto, Odd graphs and Selmer groups of certain elliptic curves, Algebra and Computation, report collection 6(2005), 10pp.
	\bibitem{new1} T.Hibino and M. Kan, $\theta$-congruent numbers and Heegner points, Arch. Math., 77(2001), 303-308.
	\bibitem{6} N. Koblitz, Introduction to Elliptic Curves and Modular Forms, GTM97, Springer, Berlin, 1984.
	\bibitem{7} M. Laczkovich, Rational points of some elliptic curves related to the tilings of equilateral triangle, https://doi.org/10.1007/s00454-019-00143-5, 2019
	\bibitem{14} J. Pan and Y. Tian, Toric periods and non-tiling numbers, preprint, 2020.
    \bibitem{8} L. R\'{e}dei, Arithemetischer Beweis des Satzes \"uber die Anzahl der durch 4 teibaren Invarianten absoluten Klassengruppe in quadratishen Zahlk\"orper, J. Reine Angew. Math. 171(1935), 55-60.
    \bibitem{9} J. H. Silverman and J. Tate, Rational Points on Elliptic Curves, Undergraduate Texts in Math., Springer, 1992.
    \bibitem{10} Y. Tian, Congruent Numbers and Heegner Points, Cambridge Jour. of Math., 2(2014), 117-161.
    \bibitem{11} Y. Tian, X. Yuan and S. Zhang, Genus periods, genus points and congruent number problem, Asian J. Math. , 21(2017), 721-774.
    \bibitem{12} J. Top and N. Yui, Congruent numbers and their variants, Algorithmic Number Theory, MSRI Publications, 44(2008), 613-639.
    \bibitem{13} S. Yoshida, Some variants of the congruent number problem, I, II, Kyushu J. Math., 55(2001), 387-404, 56(2002), 149-165.

	
	
\end{thebibliography}
\end{document}